\documentclass{amsart}
\usepackage{amssymb}
\usepackage[T1]{fontenc}
\usepackage{yfonts}

\newtheorem{thm}{Theorem}[section]
\newtheorem{cor}[thm]{Corollary}
\newtheorem{lem}[thm]{Lemma}

\theoremstyle{definition}

\theoremstyle{remark}

\numberwithin{equation}{section}
\newcommand{\R}{\mathbb R}
\newcommand{\eps}{\varepsilon}

\newcommand{\rt}{\rightarrow}

\newcommand{\ca}{\mathcal {C}}

\newcommand{\g}{\bf{g}}
\newcommand{\s}{\mathbf{s} \mathbf{o}}

\newcommand{\C}{{\mathbb C}}

\begin{document}

\title[On Wilking's Criterion for the Ricci flow]{On Wilking's criterion for the Ricci flow}
%\author{Atreyee Bhattacharya}
\author{H. A. Gururaja}
\author{Soma Maity}
\author{harish seshadri}
\address{department of mathematics,
Indian Institute of Science, Bangalore 560012, India}

\email{gururaja@math.iisc.ernet.in}
\email{harish@math.iisc.ernet.in}
\email{somamaity@math.iisc.ernet.in}

%\subjclass{53C21}

%\keywrds{Weyl Curvature, Euler Characteristic,
%Chern-Gauss-Bonnet Theorem, Asymptotically Flat Manifolds, Yamabe
%metric.}
%\thanks{This work was supported by DST Grant No. SR/S4/MS-283/05}
%\date{}%
%\dedicatory{Dedicated to Pierre B\'erard and Sylvestre Gallot on the occasion of their 60th birthdays}%
%\commby{}%
% ----------------------------------------------------------------
\begin{abstract}

 B Wilking has recently shown that one
can associate a Ricci flow invariant cone of curvature operators
$C(S)$, which are nonnegative in a suitable sense, to every
$Ad_{SO(n,\C)}$ invariant subset $S \subset {\bf so}(n,\C)$. For
curvature operators of a K\"ahler manifold of complex dimension
$n$, one considers $Ad_{GL(n,\C)}$ invariant subsets $S \subset
{\bf gl}(n,\C)$. In this article we show: \vspace{2mm}

(i) If $S$ is an  $Ad_{SO(n,\C)}$ subset, then $C(S)$ is contained in the
cone of curvature operators with nonnegative isotropic curvature
and if $S$ is an  $Ad_{GL(n,\C)}$ subset, then $C(S)$ is contained in the
cone of  K\"ahler curvature operators with nonnegative orthogonal
bisectional curvature. \vspace{2mm}

(ii) If $S \subset {\bf so}(n,\C)$ is a closed $Ad_{SO(n,\C)}$
invariant subset and $C_+(S) \subset
 C(S)$ denotes the cone of curvature operators which are {\it positive} in the appropriate sense
 %and $M$ a Riemannian manifold whose curvature operator lies in $C_+(S)$ at every point of $M$, we write  then
then one of the two  possibilities holds: (a) The connected sum of
any two Riemannian manifolds with curvature operators in $C_+(S)$
also admits a metric with curvature operator in $C_+(S)$ (b) The
normalized Ricci flow on any compact Riemannian manifold $M$ with
curvature operator in $C_+(S)$ converges to  either to a metric of
constant positive sectional curvature or constant positive
holomorphic sectional curvature or $M$ is a rank-1 symmetric
space.

\end{abstract}
\maketitle
% ----------------------------------------------------------------
\section{Introduction}
In ~\cite{wil}, B. Wilking gave a simple criterion for the positivity of certain curvatures to be preserved under the Ricci flow. Wilking's criterion captures many of the
known Ricci flow invariant positivity conditions. Some examples are positive curvature operator ~\cite{ham}, positive complex sectional curvature and positive isotropic curvature and its variants ~\cite{bs}, ~\cite{bre1}.

Wilking's criterion can be stated for general Lie algebras. Here
we confine ourselves to case of the Lie algebra being ${\bf
so}$$(n,\C)$ or ${\bf gl}$$(n,\C)$:  Let $(M,g)$ be a compact
Riemannian $n$-manifold. For any $p \in M$, choose an isometry
$L_g: T_pM \rt {\mathbb R}^n$ where $\R^n$ is endowed with the
standard metric. One then has an identification of   $\wedge^2
T_pM$ with $\g:=\s$$(n)$, the Lie algebra of skew-symmetric
endomorphisms of ${\mathbb R}^n$. This identification is via the
composition of $L_g \wedge L_g$ with the map $\phi: \wedge^2 \R^n
\rt {\bf so} (n)$ given  by
\begin{equation}\label{one}
 \phi(u \wedge v)(x)= \langle u,x \rangle v - \langle v, x \rangle u \ \ \ x \in {\mathbb R}^n.
\end{equation}

 Consider the space $S^2(\g)$ of endomorphisms of $\g$ which are skew-symmetric with
respect to the inner-product $\langle A, B \rangle = - \frac{1}{2}tr (AB)$.  Extend this inner product to a Hermitian form on $\g \otimes_\R \C = {\bf so}$$(n,\C)$ and extend any $R \in S^2(\g)$ to a complex linear map $ \bf {so}$$(n,\C) \rt  \bf {so}$$(n,\C)$. Denoting the extensions by the same symbols, one has ~\cite{wil} \vspace{2mm}

{\it Let S be a subset of the complex Lie algebra $\bf {so}$$(n,\C)$. If S is invariant under the adjoint action of
the corresponding complex Lie group $SO(n,\C)$, then the set
$$C(S)= \{ R \in S^2({\bf g}) \ \vert \ \langle R(X),  X \rangle \ge 0 \ \ {\rm for \ all} \ X \in S \}$$
is invariant under the ODE
$$ \frac {dR}{dt} = R^2 + R^\#.$$} \vspace{2mm}

Next consider the case of K\"ahler manifolds. Let $(M,J,g)$ be a K\"ahler manifold (with $J$ denoting the almost complex structure) of complex dimension $n$ and $p \in M$.
$J:T_pM \rt T_pM$ extends to an operator $\hat J$ of $\wedge^2 T_pM$ satisfying $\hat J ^2=I$ with $I$ being the identity. If $V_p$ denotes the $+1$ eigenspace of $\hat J$, then
$V_p$ is an invariant subspace of the curvature operator of $(M,g)$. Choose a complex linear isometry
$L_g: T_pM \rt \C^n=\R^{2n}$. Then the map $\phi \circ L_g \wedge L_g : \wedge^2 T_p \rt {\bf so}(2n)$  restricts to an isomorphism $V_p \rt {\bf u}(n)$ where
${\bf u}(n)$ is identified with the set of skew-symmetric endomorphisms of ${\mathbb R}^{2n}$ which commute with the standard almost complex structure of  ${\mathbb R}^{2n}$. We again
consider $R \in S^2({\bf u}(n))$ and the complex-linear extension of $R$ to ${\bf u}(n) \otimes \C = {\bf gl}(n, \C)$. For $S \subset {\bf gl}(n, \C)$ which is $Ad_{ GL(n, \C)}$ invariant, Wilking's
result asserts that the cone $C^K(S)$ of K\"ahler curvature operators nonnegative on $S$ is invariant under
$ \frac {dR}{dt} = R^2 + R^\#$. \vspace{2mm}

In the above framework, the cone of curvature operators with nonnegative isotropic curvature corresponds
is the cone $C(S_0)$ where
$$S_0 = \{ X \in { \bf so}(n,\C) \ \vert \  rank(X)=2, \ \ X^2=0 \}$$
and the cone of K\"ahler curvature operators with nonnegative orthogonal bisectional curvature is the cone
$C^K(S_1)$ where
$$S_1= \{ X \in { \bf gl}(n,\C) \ \vert \ rank(X)=1, \ \ X^2=0 \}.$$

\begin{thm}\label{mai1}
(i) For any $Ad_{SO(n, \C)}$ invariant subset $S \subset { \bf so}(n,\C)$ one has
 $$C(S) \subset C(S_0).$$ 

Moreover, if $\langle R(X),X \rangle=0$ for all $X \in S$ then 
 $\langle R(Y),Y \rangle=0$ for all $Y \in S_0$. \vspace{2mm}

(ii) If $S \subset { \bf gl}(n,\C)$ is $Ad_{GL(n, \C)}$ invariant then
$$C^K(S) \subset C^K(S_1).$$

Again, if  $\langle R(X),X \rangle=0$ for all $X \in S$ then 
 $\langle R(Y),Y \rangle=0$ for all $Y \in S_1$.
\end{thm}
\vspace{3mm}

To state our second result, consider the complex linear extension $\phi : \wedge ^2 \C^n \rt {\bf so}(n,\C)$ of the map $\phi$ in (\ref{one}). We
call an element of the form $\phi(u \wedge v) \in {\bf so}(n,\C)$ a {\it simple} element. \vspace{2mm}

Let $S$ be a closed $Ad_{SO(n, \C)}$-invariant subset of ${\bf so}(n,\C)$ and
$$C_+(S):=\{ R \in S^2({\bf g}) \ \vert \ \langle R(X),  X \rangle >  0 \ \ {\rm for \ all} \ X \in S \}.$$
If $(M,g)$ is a  Riemannian manifold whose curvature operator (at
every point of $M$) lies in $C_+(S)$, we say that  $(M,g)$ has
{\bf positive S-curvature}. \vspace{2mm}

Let ${\mathcal C}$ be the class of all closed $Ad_{SO(n, \C)}$-invariant subsets of ${\bf so}(n,\C)$
and $\ca_0$ be the subcollection consisting of those invariant sets which contain a nonzero simple element
of the form $\phi(e \wedge u)$ with $e \in \R^n$,  $u \in \C^n$ and $(e,u)=0$. Here $(\ , \ )$ denotes the complex bilinear extension of the usual inner product on $\R^n$ to $\C^n$.

\begin{thm}\label{mai2}
Let $S$ be a closed $Ad_{SO(n, \C)}$-invariant subset of ${\bf so}(n,\C)$. \vspace{2mm}

(i) If $S \in \ca \setminus \ca_0$ then the connected sum of any
two Riemannian manifolds with positive $S$-curvature also admits a
metric with positive $S$-curvature. \vspace{2mm}

(ii) If $S \in  \ca_0$ and $M$ is any compact Riemannian manifold
with positive $S$-curvature then one of the following holds:

(1) The normalized Ricci flow on $M$ converges to a metric of constant positive sectional curvature.

(2) $M$ is K\"ahler and the normalized Ricci flow on $M$ converges to a metric of constant positive holomorphic  sectional curvature.

(3) $M$ is isometric to a rank-1 symmetric space.
\end{thm}
\vspace{2mm}

It should be noted that this result is sharp in the sense that for
any $S \in \ca \setminus \ca_0$ there is a Riemannian manifold
$(M,g)$ with positive $S$-curvature such that the normalized Ricci
of $(M,g)$ dose not converge to a metric on $M$. Indeed,  the
standard product manifold $S^{n-1} \times S^1$ has positive
$S$-curvature for any $S \in \ca \setminus \ca_0$. \vspace{2mm}

Moreover for any  $S \in  \ca_0$ there are compact Riemannian manifolds $M_1$ and $M_2$ with positive $S$-curvature whose connected  sum does not admit a metric with positive $S$-curvature: we can
take $M_1=M_2$ be a non-simply connected quotient of $S^n$ with canonical metric. Since
positive $S$-curvature implies positive Ricci curvature when  $S \in  \ca_0$ (see the proof
of Corollary \ref{ddd}), the Myers-Bonnet theorem on the finiteness of fundamental group rules
out the existence of a metric with positive $S$-curvature on $M_1 \# M_2$. \vspace{2mm}

A corollary of the proof of Theorem \ref{mai2} is the following:

\begin{cor}
Let $S$ be an $Ad_{SO(n, \C)}$-invariant subset of ${\bf so}(n,\C)$. 

The normalized Ricci flow $g(t)$ of any compact Riemannian manifold $M$ with positive $S$-curvature
converges to a metric of constant positive sectional curvature if and only if  the product metric
$g(t)+ds^2$ on $M \times \R$ has positive isotropic curvature for any $t>0$.
\end{cor}

The ``if'' part of the above theorem is a result of S. Brendle ~\cite{bre1}. \vspace{2mm}

{\bf Note:} After this article was written  we came to know from B. Wilking that he was aware of Theorem \ref{mai1}, even though it was not mentioned in ~\cite{wil}.

\section{Minimal Ricci flow invariant positivity conditions}
\noindent We prove Theorem \ref{mai1} in this section. \vspace{2mm}

\noindent {\bf Proof} (of Part (i) of Theorem \ref{mai1}):

\begin{lem}\label{qw}
Let $X \in S$. If $X^2 \neq 0$ then there exists a sequence $\{P_1,P_2,..\} \subset SO(n,\C)$ such that
$$ lim_{k \rt \infty} \Vert P_k X P_k^{-1} \Vert = \infty .$$
\end{lem}

\noindent {\bf Proof:}
Let $(\ , \ )$ denote the complex linear extension of the standard inner product on $\R^n$ to $\C^n$. An endomorphism  $P$ of $\C^n$ belongs to
$SO(n,\C)$ if and only if it maps  orthonormal bases with respect to $( \ , \ )$ to each other. Since $X^2\neq 0$, there ia a $v \in \C^n$ such that $Xv$ is not isotropic i.e. $(Xv,Xv) \neq 0$. By perturbing $v$, if necessary, we can assume that $v$ is not isotropic as well.
Let $L:=span_\C \{ v \}$ and $L^\perp$ the orthogonal complement of $L$ with respect to $( \ , \ )$.
Since $X \in {\bf so}(n,\C)$, we have $(v,\ Xv) =0$ i.e. $Xv \in L^\perp$. We can find $\{ v_3,...,v_n\}$ such that
$\{ c Xv , \ v_3,...,v_n\}$ is an orthonormal basis of $L^\perp$ with respect to $( \ , \ )$ with
$c= (Xv, \ Xv)^{- \frac{1}{2}}$.  Choose another $( \ , \ )$ orthonormal basis  $\{x_k,w_3,...,w_n\}$ of $L^\perp$ such that the Euclidean norm $\Vert x_k \Vert = k$. Define $P_k \in SO(n,\C)$ by
$$ P_k(v)=v, \ \ P_k(c Xv) = x_k, \ \ P_k(v_i)=w_i, \ \ i=3,..,n.$$
Then $\Vert P_k X P_k^{-1}v \Vert =c^{-1}  \Vert x_k \Vert \rt \infty$ as $k \rt \infty$. Hence
$\Vert P_k X P_k^{-1} \Vert  \rt \infty$. \hfill $\square$ \vspace{2mm}

Now we continue with the proof of (i), Theorem \ref{mai1}. Let $X
\in S$ and $X^2 \neq 0$. Let $P_k$ be as in  Lemma \ref{qw} and
$\lambda_k := \Vert P_k X P_k^{-1} \Vert^{-1}$. Then a subsequence
of $T_k = \lambda_k P_k X P_k^{-1}$ converges to a nonzero $T \in
{\bf so}(n,\C)$. It can be readily checked that if $p$ is the
degree of the minimal polynomial of $X$, then $T^p=0$.

Note that $\langle R(T),T \rangle = \vert \lambda_k \vert^2 lim_{k \rt \infty} \langle R( P_k X P_k^{-1}), \ P_k X P_k^{-1} \rangle \ge 0$. In fact, if $T_1 \in O(T)$, where $O(T)$ denotes the adjoint orbit of $T$, then
$\langle R(T_1),T_1 \rangle \ge 0$.

The proof is completed by the following classical facts (cf ~\cite{nil}, Theorems 4.3.3 and Proposition 5.4.1): Let ${\bf g}$ be a  simple Lie algebra. Then there exists a nonzero nilpotent orbit of minimal dimension, denoted by $O_{min}$, which is contained in the closure of any
nonzero nilpotent orbit. For $ {\bf so}(n,\C)$, \ $O_{min}=O(A)$ for any $A$ satisfying $A^2=0$ and $rank(A)=2$ i.e. $O_{min}=S_0$.

Hence $S_0 \subset \overline {O(T)}$ and $\langle R(v), v \rangle \ge 0$ for all $v \in S_0$ by continuity.

If $X^2=0$ to begin with, then again $S_0 \subset \overline {O(X)}$ and continuity gives the desired conclusion.
\hfill $\square$ \vspace{3mm}

{\bf Proof} (of part (ii) of Theorem \ref{mai1}): Let $X \in S$. As in Lemma \ref{qw} we can find $P_k \in GL(n, \C)$ such that $lim_{k \rt \infty} \Vert P_k X P_k^{-1} \Vert = \infty$.
As before let $\lambda_k^{-1}=  \Vert P_k X P_k^{-1} \Vert$. A subsequence of $\lambda_k \Vert P_k X P_k^{-1} \Vert$  converges to a nonzero nilpotent
$Y \in {\bf gl}(n, \C)$.

Let $$J_k = \begin{pmatrix}
  0 & 1 &  &     \\
    & 0 & \ddots &     \\
    &   & \ddots & 1     \\
    &   &        &        0
 \end{pmatrix}$$ denote an elementary nilpotent Jordan block of size $k$. Then the Jordan canonical forms of $Y$ and
 any  $Z \in S_1$ can be written as
\begin{equation}\notag
Z =P^{-1} \left( \begin{array}{cc}
J_2 & 0  \\
0 & {\bf 0}   \end{array} \right) P  \ \ \ \ \ \ {\rm and} \ \ \ \ \ \ Y=Q^{-1}
 \begin{pmatrix}
  J_{k_1} &  &  &  \\
          & J_{k_2} &  &  \\
         &   & \ddots &   \\
   &  &  & J_{k_p}
 \end{pmatrix} Q
\end{equation}
where $k_1 \ge ... \ge k_p$.

Since each $J_k$ is similar to any multiple of itself, we can find $Q_k \in GL(n, \C)$ such that $Q_k^{-1}YQ_k \rt Y_1$ as $k \rt \infty$, where
$$ Y_1 =   \left( \begin{array}{cc}
J_{k_1} & 0  \\
0 & {\bf 0}   \end{array} \right).$$

Moreover, if $a_1,...,a_{k_1} \in \C^\ast$, then $J_{k_1}$ is similar to
 $$\begin{pmatrix}
  0 & a_1 &  &     \\
    & 0 & \ddots &     \\
    &   & \ddots & a_{k_1}     \\
    &   &        &        0
 \end{pmatrix}.$$
Letting $a_1 \rt 1$ and $a_i \rt 0$ for $i=2,...,k$, we see that  $\langle R(Z),Z \rangle \ge 0$ for any $R \in C(S)$. \hfill $\square$

\section{Positive curvature on connected sums}
In this section we prove (i) of Theorem \ref{mai2}.

Let $S$ be an $Ad_{SO(n,\C)}$ invariant subset of ${\bf so}(n,\C)$. We
give a criterion for $S$ which will ensure that if $M_1$ and $M_2$
are Riemannian manifolds with positive $S$-curvature then the
connected sum $M_1 \# M_2$ admits a Riemannian metric with
positive $S$-curvature.

A natural procedure for constructing such a metric is to remove
points $p_i \in M_i$ and deform the metric $g_i$  to new metrics
$\tilde g_i$ on $M_i \setminus \{ p_i \}$ such that (i) $\tilde
g_i$ has positive $S$-curvature and (ii) $\tilde g_i$ equals the
standard product metric on a cylinder near $p$. One can then
connect the two cylinders (possibly after rescaling $g_i$) to get
the required metric on $M_1 \# M_2$.  Note that a necessary
condition on $S$ for this to work is that the cylinder with the standard product metric should have positive $S$-curvature.

Let $e_0$ and $\{e_1,...,e_{n-1}\}$ denote orthornormal vectors tangent (at some point) to the first
and second factors of the cylinder $\R \times S^{n-1}$. It is easily checked that $\langle R_0(X),X \rangle =0$
iff $X=c_1 X(e_0,e_1)+...+c_{n-1}X_(e_0,e_{n-1})$ for some $c_i \in \C$.  Here we use the following notation:
\begin{equation}\label{two}
{\rm If} \ u,v \in \C^n \ {\rm then} \ X(u,v):= \phi(u \wedge v)
\end{equation}
where $\phi$ is the map defined  in (\ref{one}). If a real orthonormal basis $\{e_0,...,e_{n-1}\}$ has been fixed we put
$$X_{ij}:=X(e_i,e_j).$$

The necessary condition on $S$ in the previous paragraph is that
it does not contain any element of the form $X=c_1 X(e_0,e_1)+...+c_{n-1}X_(e_0,e_{n-1})$  . The main result
of this section is that this condition is sufficient. We first
note that an $Ad_{SO(n,\C)}$ invariant subset $S \subset {\bf
so}(n,\C)$ contains a nonzero element of the form  $c_1
X(e_0,e_1)+...+c_{n-1}X(e_0,e_{n-1})$ for some $c_i \in \C$  and
some real orthonormal basis $\{e_0,...,e_{n-1}\}$ if and only if it
contains a nonzero simple element of the form $X(e,
 u)$ with $e \in \R^n$,  $u \in \C^n$ and $(e,u)=0$ .

\begin{thm}\label{pol}
Let $S$ be an $Ad_{SO(n,\C)}$ invariant subset of $ {\bf
so}(n,\C)$ such that the closure $\overline S$  does not contain
any nonzero simple element of the form $X(e, u)$ with $e \in
\R^n$, $u \in \C^n$ and $(e,u)=0$.

 If $(M_i,g_i)$, with $i=1,2$, are
compact Riemannian manifolds with positive $S$-curvature then the
connected sum $M_1 \# M_2$ admits a Riemannian metric  with
positive $S$-curvature.
\end{thm}

%{\bf Remark:} Note that we do not need $Ad_{SO(n,\C)}$ invariance of $S$ for the above result. If one %does have this
%invariance then the hypothesis on $S$ is equivalent to requiring that $S$ does not contain an element of %the above form
%for {\it some} orthonormal basis of $\R^n$.
%\begin{lem}
%Let $S$ be an $Ad_{SO(n,\C)}$ invariant subset of $ so(n,\C)$ and $(M,g)$ be a compact Riemannian manifold with curvature operator in $C_+(S)$.
%\end{lem}

{\bf Proof:} The construction of the metric proceeds along the
same lines as the Micallef-Wang ~\cite{mw} construction of positive isotropic
curvature on connected sums.  Let $(M,g)$ be a Riemannian manifold
with positive $S$-curvature. Let $p \in M$ and $B$ be a geodesic
ball centered at $p$ and radius less than the injectivity radius
at $p$. We will construct a smooth positive function $u$ on $M
\setminus \{p\}$ such that $u\vert_{M \setminus B} \equiv 1$ and
near $r=0$, $u^2g$ is close to the product metric on $\R^+ \times
S^{n-1}(\rho)$ in the $C^2$ topology, for a suitable choice of
$\rho$. This can be done so that $u^2g$ has positive
$S$-curvature.
 One can then deform this to the product metric again maintaining positive $S$-curvature.

 Let $\{e_1,...,e_n\}$ be an orthonormal frame for $g$. If $\tilde g =u^2 g$ and $\tilde e_i= \frac {1}{u} e_i$ then
$ \{\tilde e_1,...,\tilde e_n\}$ is an orthornomal frame for  $\tilde g$. One then has
\begin{align}
\tilde R_{ijkl} :=& \ \ \tilde g (\tilde R (\tilde e_i \wedge \tilde e_j),
\tilde e_k \wedge \tilde e_l) \\
=& \ \ u^{-2} g(R(e_i \wedge e_j), e_k \wedge e_l) \ \ + \\
& u^{-2}\Bigl (g \odot (u^{-1}Hess(u) -2u^{-2}du \otimes du +
\frac{1}{2} u^{-2}\vert du \vert^2g ) \Bigr )(e_i,e_j,e_k,e_l),
\notag
\end{align}
where  $\odot$ denotes the Kulkarni-Nomizu product of symmetric 2-tensors.

If $u=u(r)$ is a radial function, then
$$ Hess(u) = u^{\prime \prime} dr \otimes dr - u^{\prime} h$$
where $h$ is the extension of the second fundamental form of the geodesic spheres $S_r$ to $B$ obtained by
putting $h(v, \frac {\partial}{\partial r})=0$ for any tangent vector $v$.

Note that if $\{e_1,...,e_{n-1}\}$ is an orthonormal basis of
$T_qS_r$ which diagonalizes the second fundamental form of $S_r$
at $q$ and $e_0 = \frac {\partial}{\partial r}$, then

\begin{equation}\label{as}
\tilde R_{ijkl} = \frac {1}{u^2} R_{ijkl}
\end{equation}
whenever $\{i,j\} \neq \{k,l\}$ with $ 0 \ \le \ i, j,k,l \ \le \
n-1$. \vspace{2mm}

%(1) $i=0, \ \ 1 \ \le \ j,k,l \ \le \ n-1$ \ \ or \ \ $i=l=0, \ \ j \neq k$

%(2) $ 1 \ \le \ i,j,k,l \ \le \ n-1$ \ and \ $\{i,j\} \neq \{k,l\}$. \vspace{2mm}

Also, if we let $K_{ij} =R(e_i,e_j,e_j,e_i)$, etc then
\begin{equation}
\tilde K_{ij} = u^{-2} \Bigl ( K_{ij} - \Bigl (\frac {u^\prime}{u} \Bigr)^2 - \Bigl ( \frac {u^{\prime \prime}}{u} - 2 \Bigl (\frac {u^\prime}{u} \Bigr )^2  \Bigr ) (e_i(r)^2 + e_j(r)^2 )
+ \Bigl (\frac {u^\prime}{u}\Bigr )(h_{ii}+h_{jj}) \Bigr )
\end{equation} \vspace{2mm}

Noting that the $h_{ii}$ are principal curvatures of $S_r$, we have

\begin{equation}\label{a}
 h_{ii} = - \frac{1}{r} + O(r).
\end{equation}

Put
\begin{equation}\label{b}
\frac {u^\prime}{u} = - \frac {\alpha(r)}{r}
\end{equation}
where $\alpha(r) \equiv 0$ near $r=0$, \ $\alpha(r) \equiv 1$ near $r=1$ and $\alpha' \le 0$. \vspace{2mm}

It follows from (\ref{a}) and (\ref{b}) that there exists $C>0$ such that
\begin{equation}\label{h}
\tilde K_{ij} \ \ge \ u^{-2} \Bigl ( K_{ij} + \frac {\alpha (2-\alpha)}{r^2} -C \alpha  \Bigr ) \ \ \ \ {\rm for} \ \ \ \ i,j \ge 1
\end{equation}
and
\begin{equation}\label{i}
\tilde K_{0j} \ \ge  \ u^{-2} \Bigl ( K_{0j} + \frac {\alpha ^{\prime}}{r} -C \alpha  \Bigr ) \ \ \ \ {\rm for} \ \ \ \ j \ge 1 \vspace{5mm}
\end{equation}

Now let $X = \sum_{0 \le i <j} a_{ij} X_{ij}$ with $a_{ij} \in \C$. We write
$$X=X_1 +X_2$$
with $X_1 = \sum_{1 \le i < j } a_{ij}X_{ij}$ and $X_2=X-X_1$. It
follows from the hypothesis of Theorem ~\ref{pol} that
\begin{align}\label{lkj}
 & there \ \ exists \ \ k >0  \ \ such  \ \ that \ \ if \ \ X \in S \ \
 and \ \ \Vert X \Vert^2 = \sum_{0 \le i < j } \vert a_{ij} \vert^2 =1 \\ \notag
 & then \ \ \Vert X_1 \Vert \ge k.
\end{align}

Recalling the notation in Page 1, we have, by definition
$$ \langle R(X),X \rangle \ := \  \langle R \circ \rho_g^{-1} (X), \ \rho_g^{-1}(X) \rangle_g,$$
where 
\begin{equation}\label{nnn}
\rho_g = \phi \circ L_g \wedge L_g
\end{equation} 
Hence, if $X \in S$
and $\Vert X \Vert =1$

\begin{align}\notag
 \langle \tilde R(X), X \rangle \ =& \ \sum_{0 \le i<j, \ 0 \le k < l} a_{ij} \bar a_{kl} \tilde R_{ijkl}
 \\ \notag
                         \ =& \ \sum_{0 \le i <j} \vert a_{ij} \vert ^2 \tilde K_{ij} + \sum_I a_{ij} \bar a_{kl} \tilde
                         R_{ijkl} \notag
\end{align}
where $I$ is the index set  $I= \{ (i,j,k,l): \  0 \le i < j, \ 0
\le k < l \ \ {\rm and} \ \ \{i,j\} \neq \{k,l\} \}$. It follows
from \eqref{as}, \eqref{h}, \eqref{i} and \eqref{lkj} that

\begin{align}
 \langle \tilde R(X), X \rangle \ \ \ge & \ \ u^{-2} \langle  R(X), X \rangle \ +
                         u^{-2} \Bigl ( \frac {\alpha'}{r} + \frac {k \alpha (2-\alpha)}{r^2} -D \alpha \Bigr )
\end{align}
for some $D>0$.

Let $0 < \eps < 1$ \ and \  $ r_0= \sqrt \frac{k(1-\eps)}{D}$ \ so that $ 0 \le r \le r_0$ implies
\begin{equation}
\frac {k \alpha (2-\alpha)}{r^2} -D \alpha \ \ge \  \frac {k \alpha (1+\eps-\alpha)}{r^2}.
\end{equation}
Let $k_0>0$ be such that
\begin{equation}\label{ljh}
  \langle  R_p(X), X \rangle \ge k_0 \ \ {\rm for \ all } \ \ \vert X \vert =1 \ \ {\rm and} \ \ p \in M,
\end{equation}
where $R_p$ denotes the curvature operator at $p \in M$.

Combining, we have
\begin{equation}
\langle \tilde R(X), X \rangle \ \ \ge  \ \
                         u^{-2} \Bigl ( k_0 + \frac {\alpha'}{r} + \frac {k \alpha (1+\eps-\alpha)}{r^2}  \Bigr )
\end{equation}
Hence, in order to maintain positivity of $\langle \tilde R(X), X \rangle$ it is enough to find a nonincreasing
smooth function $\alpha$ on $[0, \ 1]$ which satisfies

$$ k_0 + \frac {\alpha'}{r} + \frac {k \alpha (1+\eps-\alpha)}{r^2} >0, \ \ \ \ \ \  \alpha \equiv 1 \ \ {\rm near} \ \ 0, \ \ \ \ \ \alpha \equiv 0 \ \ {\rm on} \ \ [r_0, \ 1].$$\vspace{1mm}

If one makes the change of variables $r=r_0e^{-t}$ then finding such an $\alpha$ is equivalent to finding  a
smooth nondecreasing $\beta:[0,\infty) \rt [0,\infty)$ such that
$$ \beta' < k_0r_0e^{-2t} + k \beta(1+\eps -\beta), \ \ \ \ \ \ \beta \equiv 0 \ \ {\rm near} \ \ 0, \ \ \ \ \ \beta \equiv 1 \ \ {\rm for \ sufficiently \ large } \ t.$$
One can construct such a $\beta$ by modifying near $0$ and $\infty$ the solution of $\beta' = k \beta(1+\eps -\beta)$ which has the form $$\beta(t)=\frac {(1+\eps)Ae^{(1+\eps)kt}}{1+Ae^{(1+\eps)kt}}$$
for a sufficiently small constant $A$.

With this choice of $\beta$ and consequently $u$, the metric $\tilde g = u^2g$ is a complete metric on
$M -\{ p \}$ which has positive $S$-curvature, is the same as $g$ outside $B$, and becomes $C^2$ close to the
product metric $dt^2 + \rho^2 g_0$ on $\R^+ \times S^{n-1}$ for some $\rho >0$.  Finally, we perform the above deformations on the manifolds $M_i$ and points $p_i \in M_i$ to get complete metrics with
positive $S$-curvature on $M_i - \{ p_i \}$ which are arbitrarily close to the product metric on
$dt^2 + \rho^2 g_0$ on $B_i - \{ p_i \} = \R^+ \times S^{n-1}$. One can then connect the two cylinders smoothly maintaining positive $S$-curvature. \hfill $\square$

\section{Convergence of Ricci flow}
\noindent We prove (ii) of Theorem \ref{mai2}. \vspace{2mm}

 %Let $V$ be an n dimensional inner product space over $\mathbb{R}$ and $\{e_i\}$
%an orthonormal basis. Then identify $\Lambda^2V$ to $so(n)$ by mapping $e_i\Lambda e_j$ to the matrix $E_{ij}$, %where ij th and ji th entries of $E_{ij}$ are $1$ and $-1$ and other entries are $0$.Next we consider %$V\otimes\mathbb{C}$ and identify $\Lambda^2V\otimes\mathbb{C}$
%with $so(n,\mathbb{C})$.

%Let $S_2$ be the $Ad_{SO(n,\mathbb{C})}$ invariant subset of ${\bf so}(n,\mathbb{C})$ given by
%$$S_2 = \{ X \in SO(n,\mathbb{C}) \ \vert \  Trace(X^3)=0, \ Rank(X)=2 \}.$$
%Positivity of the curvature operator of $M$ on $S_2$ is equivalent to positive isotropic curvature  on $M \times \R$.

As usual, $S$ will denote an arbitrary $Ad_{SO(n,\mathbb{C})}$
subset of ${\bf so}(n,\mathbb{C})$. We can assume that $S$ is
scale-invariant, i.e. $S = \C^\ast S:= \{ \lambda X \ : \ X \in S,
\ \lambda \in \C^\ast \}$, since a curvature operator $R \in
C_+(S)$ iff $R \in C_+(\C^\ast S)$.

%\subsection{A canonical form for elements of $S_2$}.
%\Bea C(S):=\{R\in S^2(so(n,\mathbb{C})|R(v,\overline{v})\geq 0, \forall v \in S \}
%\Eea Let $(M,g)$ be a riemannian manifold. We can identify $T_pM$
%with $V$ and then the curvature operater $R\in S^2(so(n))$. Extend $R$
%complex bilinearly.

%\begin{thm}Let $M$ be a closed manifold which admits a riemannian metric $g$
%such that at every point the curvature operator belongs to $C(S)$.Then $M$ is
%diffeomorfic to sphere.
%\end{thm}

\begin{lem}\label{pcb}
If $S$ contains a nonzero simple element of the form $X(e, u)$ with $e \in \R^n$,  $u \in \C^n$ and $(e,u)=0$ then one of
the following possibilities holds:

(i) $S$ contains all elements of the form $X(v_1,v_2)$ where
$v_1,v_2$ are $( \ , \ )$ orthonormal vectors in $\C^n$.

(ii) $S$ contains all elements of the form $X(v,e_1)+iX(v,e_2)$
where $v \in \C^n$, $e_1, e_2 \in \R^n$ and $v,e_1,e_2$ are $( \ ,
\ )$ orthonormal.
\end{lem}

\noindent {\bf Proof:}  Let $X(e,u) \in S$ with  $e \in \R^n$,  $u \in \C^n$ and $(e,u)=0$.
Since $S$ is scale invariant, we can assume that $(e,e)=1$.

Suppose $u$ is not isotropic, i.e. $(u,u) \neq 0$. Then we can assume that $(u,u)=1$.
Given $( \ , \ )$ orthonormal $v_1,v_2 \in \C^n$
we can then find $P \in SO(n,\C)$ with $P(e)= v_1, \ P(u) =v_2$.
Since $P X(e,u)P^{-1} = X(Pe,Pu)=X(v_1,v_2)$ we see that (i) holds  if $u$ is not isotropic.

Suppose $u$ is isotropic. Then there exist  orthonormal $f_1, \
f_2 \in \R^n$ such that  $u=f_1+ if_2$. Note that $\{ e,f_1,f_2
\}$ is an orthonormal set since $(e,u)=0$ by hypothesis. We have
$X(e,u)= X(e,f_1)+iX(e,f_2)$. Let $v \in \C^n$, $e_1, e_2 \in
\R^n$ with $v,e_1,e_2$ are $( \ , \ )$ orthonormal.

Let $P \in SO(n,\C)$ be an element satisfying
$$P(e)=v, \  P(f_1)=e_1, \ P(f_2)=e_2.$$
Then
$$ P^{-1} XP = X(v,e_1)+iX(v,e_2).$$

This proves (ii) of Lemma \ref{pcb}. \hfill $\square$
\vspace{3mm}

\begin{cor}\label{ddd}
Let $M$ be a Riemannian manifold with positive $S$-curvature where
$S$ is an $Ad_{SO(n,\C)}$ invariant subset   containing  a nonzero
simple element of the form $X(e, u)$ with $e \in \R^n$, $u \in
\C^n$ and $(e,u)=0$. $M \times \R$ then has nonnegative isotropic
curvature. Moreover $M$ is locally de Rham irreducible.

%Moreover if $M$ has positive isotropic
%and positive $S$-curvature then $M \times \R$ has positive isotropic curvature.

\end{cor}
\noindent {\bf Proof:} By \cite{bre1}, $M \times \R$ has
nonnegative isotropic curvature if and only if
\begin{equation}\label{eq}
K_{13} + \lambda^2 K_{14} + K_{23}+ \lambda^2K_{24} -2\lambda
R_{1234} \ge 0
\end{equation}
for all orthonormal $4$-tuples $\{e_1,e_2,e_3,e_4\}$ in $T_pM$,
for any $p \in M$ and for all $\lambda \in [-1,1]$.

Fix an orthonormal  $4$-tuple $\{e_1,e_2,e_3,e_4\}$. By Lemma
\ref{pcb} one of the following possibilities holds: \vspace{2mm}

(1) $\langle R(X(v_1,v_2)), X(v_1,v_2) \rangle  > 0$ for all $( \
, \ )$ orthonormal vectors $\{v_1,v_2\}$. For any  $\lambda , \mu
\in \R \setminus \{\pm 1\}$ let  $v_1 =
(1+\lambda^2)^{-\frac{1}{2}}(e_1 +i \lambda e_2), \
v_2=(1+\mu^2)^{-\frac{1}{2}}(e_3 + i \mu e_4)$. We get
\begin{equation}\label{ww}
 K_{13} + \mu^2 K_{14} + \lambda^2 K_{23}+ \lambda^2 \mu^2 K_{24} -2\lambda \mu R_{1234} > 0.
\end{equation}
 Letting $\mu \rt 1$ we see that (\ref{eq}) holds. \vspace{2mm}

(2) $ \langle R(X(v,e_1)+iX(v,e_2)), X(v,e_1)+iX(v,e_2) \rangle
>0$ where $v \in \C^n$ and $v,e_1,e_2$ are $( \ ,
\ )$ orthonormal. For any $\lambda \neq \pm 1$ let $v=
(1+\lambda^2)^{-\frac{1}{2}}(e_3 + i \lambda e_4)$. We then have
\begin{equation}\label{ww2}
K_{13} + \lambda^2 K_{14} + K_{23}+ \lambda^2K_{24} -2\lambda
R_{1234} > 0.
\end{equation}
Hence (\ref{eq}) follows.

The local irreducibility of $M$ is clear: Suppose that an open
subset $U \subset M$ is a Riemannian product $U=U_1 \times U_2$.
Take $e_3 \in U_1, \ e_1,e_2,e_4 \in U_2$. In case (1), putting
$\lambda=\mu=0$, we have a contradiction to (\ref{ww}). In case
(2) we get a contradiction to (\ref{ww2}) by taking $\lambda=0$.
\hfill $\square$

\vspace{3mm}

We need the following consequence, due to Wilking (Appendix of
~\cite{wil}), of the Strong Maximum Principle of Brendle-Schoen
~\cite{bs2}.

\begin{lem}\label{ooo}
Let $S$ be an $Ad_{SO(n,\C)}$ invariant subset of ${\bf so}(n,\C)$
and $(M,g)$ a compact locally irreducible Riemannian manifold with
nonnegative $S$-curvature. Let $g(t)$ be the solution to Ricci
flow starting at $g$. For any $t>0$ one of the following holds:

(i) $(M,g(t))$ has positive $S$-curvature or

(ii) $(M,g)$ is K\"ahler or

(iii) $(M,g)$ is isometric to a symmetric space.
\end{lem}

{\bf Proof:} Let $t>0$, $p \in M$  and $S_t(p) \subset \wedge^2T_p(M)$ be the subset 
corresponding to $S$ at time $t$ (in the notation of (\ref{nnn}) $S_t(p) = \rho_{g(t)}^{-1}(S)$).
Let 
$$T(t)= \bigcup_{p \in M}  \{ x \in S_t(p) \ : \ \langle R(t)(x), x \rangle _t =0 \}.$$
We claim that if the holonomy group $Hol(p)=SO(n, \R)$ then it follows that
$(M,g(t))$ has positive $S$-curvature: Otherwise there is a
nonzero $x \in T(t)$. By the Appendix of ~\cite{wil} $T(t)$ is invariant under the
action of $Hol(p)$. This translates to saying that the curvature vanishes on the
orbit of $X=\rho_{g(t)}(x)$. By (i) of Theorem \ref{mai1}, it follows that all
isotropic curvatures vanish as well. In particular the scalar
curvature of $(M,g(t))$ has to be zero. By the evolution equation
for scalar curvature 
$$\frac {\partial s}{\partial t} = \Delta s +
2 \vert Ric \vert^2$$  
it follows that the scalar curvature and
Ricci curvature of $g(t')$ have to vanish for all $t' \in [0,t]$.
By the vanishing of $s$ on $[0,t]$ it follows that the isotropic
curvature of $g(s)$ is zero and hence each $g(s)$ is  conformally
flat. Since we already know that $g(s)$ is Ricci flat, we conclude
that each $g(s)$ is flat. In particular $g=g(0)$ is flat, which
contradicts the local irreducibility of $g$. \vspace{2mm}

Moreover if $g(t)$ has positive $S$-curvature for some $t$ then
$g(t')$ has positive $S$-curvature for all $t' \ge t$. Hence we
can assume that $g(t')$ does not have positive $S$-curvature for
all $t' \in [0,t)$ for some $t>0$. By the previous paragraph
$Hol(p)$ is a proper subgroup of $SO(n)$ for all $t' \in [0,t)$.
Berger's classification of holonomy groups then leads to the
following possibilities: Either $(M,g(t))$ is a locally symmetric
space of rank $\ge 2$ or $Hol(p)=U(n)$ or $Sp(n)Sp(1)$. Hence one
can conclude the same about $g$. If $Hol(p)=Sp(n)Sp(1)$ then
$(M,g)$ is quaternion-K\"ahler and has nonnegative isotropic
curvature (by (i) of Theorem \ref{mai1}). By a result of Brendle ~\cite{bre2}
such a $(M,g)$ is locally symmetric.

If $Hol(p)=U(n)$ then $(M,g)$ is K\"ahler. \hfill $\square$
\vspace{2mm}

We can now complete the proof of Theorem \ref{mai2}. \vspace{2mm}

{\bf Proof} (of (ii), Theorem \ref{mai2}): Let $M$ be a compact
Riemannian manifold with positive $S$-curvature where $S$ contains
a nonzero simple element of the form $X(e, u)$ with $e \in \R^n$,
$u \in \C^n$ and $(e,u)=0$. By Corollary \ref{ddd}, $M$ is
irreducible and $M \times \R$ has nonnegative isotropic curvature.
Note that the latter condition is precisely that of nonnegative
$S'$-curvature where
$$S'=\{X \in {\bf so}(n,\C) \ : \ X^3=0, \ Rank(X)=2 \}.$$
If $g(t)$ is the solution to Ricci flow starting at $g$ Lemma
\ref{ooo} implies that either $g(t)$ has positive $S'$-curvature
or $(M,g)$ is K\"ahler or $(M,g)$ is locally symmetric. We now
consider the three  cases separately:

(i) By a result of Brendle ~\cite{bre1} the normalized Ricci flow starting at a
metric of  positive $S'$-curvature (i.e. with positive isotropic
curvature on $M \times \R$) converges to a metric of constant
positive sectional curvature.

(ii) Suppose $(M,g)$ is K\"ahler. We look at the two cases in the
proof of Lemma \ref{ddd}: In case (1), $M$ has positive sectional
curvature while in case (2), $K_{12}+K_{13}>0$ for any orthonormal
triple $\{e_1,e_2,e_3\}$. In either case $M$ has positive
orthogonal bisectional curvature. It  is known, by the work of X.
Chen ~\cite{che},  H. Gu - Z. Zhang ~\cite{gz} and X. Chen - S.
Sun - G. Tian   ~\cite{cst}, that the normalized Ricci flow
starting at such a metric converges to a K\"ahler metric of
constant holomorphic sectional curvature.

(iii) Let $(M,g)$ be a symmetric space with positive
$S$-curvature. Note that we have $K_{12}+K_{13}>0$ for all
orthonormal triples $\{e_1,e_2,e_3\}$ on $M$. By the work of X. S.
Liu ~\cite{liu} and N. Lee ~\cite{lee}, none of the higher rank
(i.e. rank $\ge 2$) symmetric spaces satisfy this condition. 

\hfill $\square$ \vspace{2mm}

{\bf Remark:} Cases (ii) and (iii) above actually occur: Let
$X=X(e,u)$ where $e$ is a unit vector in $ \R^n$, $u$ is isotropic
and $(e,u)=0$. Let $S$ the orbit of $X$. We claim that {\it if $M$
has weakly quarter-pinched sectional curvatures then $M$ has
$S$-positive curvature}. First note that $PXP^{-1} =X(P(e),P(u))$
where $P(e)$ is non-isotropic, $P(u)$ is isotropic and
$(P(e),P(u))=0$. Let $V \subset \C^n$ be the two-dimensional
subspace spanned by $P(e)$ and $P(u)$. We claim that there is an
orthonormal 4-frame $\{e_1,e_2,e_3,e_4\} \subset \R^n$ and $c \in
\R \setminus \{1\}$ such that
$$\{ e_1+ice_2, \   e_3+ie_4 \}$$
is a basis for $V$. To see this, one proceeds as follows: Since
$P(u)$ is isotropic, there exist orthonormal $e_3,e_4 \in \R^n$
and $a \in \R$ such that $P(u)=a(e_3+ie_4)$.  Let $v=P(e)-\alpha
P(u)$ where $\alpha \in \C$. It can be checked that one can choose
$\alpha$ so that $Re(v)$ and $Im(v)$ are perpendicular to $e_3$
and $e_4$. Since $P(u)$ is isotropic and $(P(e),P(e))=1$, $Re(v)$
and $Im(v)$ are orthogonal to each other as well.

Hence there is  an orthonormal 4-frame $\{e_1,e_2,e_3,e_4\}$ and
$c \in \R \setminus \{1\}$ such that
$$X(P(e),P(u)) = bX(e_1+ice_2, \ e_3+ie_4)$$
for some $b \in \C$. We then have
\begin{align}\notag
\vert b \vert^{-2} \langle R(PXP^{-1}), PXP^{-1} \rangle \ &= \
K_{13}+K_{14}+c^2K_{23}+c^2K_{24} -2cR_{1234}  \notag  \\
& \ \ge \frac {1}{2}(1+c^2 -2c) > 0 \notag
\end{align}
for all $c \in \R \setminus \{1\}$. We have used Berger's estimate
$\vert R_{1234} \vert \le \frac {1}{2}$ (with the normalization
$\frac{1}{4} \le K \le 1$) in the first inequality above
~\cite{ber}.

In particular, all rank-1 symmetric spaces have positive
$S$-curvature for such an $S$.

%To see this consider the subspace $W=V+\overline V$. If $V$ does
%not contain a vector in $\R^n$, then $dim_\C W=4$. Since $P(u)$ is
%isotropic we can write $P(u)=b(e_3+ie_4)$ for some orthonormal
%pair $e_3,e_4$ in $\R^n$ and $b \in \C$. We can find $e_1,e_2
%\R^n$ such that $\{e_1,e_2,e_3,e_4\}$ is an orthonormal basis of
%$W$. Since $(P(e),P(u))=0$

\end{document}